\documentclass{article}
\usepackage{color}
\usepackage{amsmath}
\catcode`\@=11 \@addtoreset{equation}{section}

\catcode`\@=12

\newtheorem{Theorem}{Theorem}[section]
\newtheorem{Proposition}{Proposition}[section]
\newtheorem{Lemma}{Lemma}[section]
\newtheorem{Corollary}{Corollary}[section]

\newcommand{\bTheorem}[1]{
\begin{Theorem} \label{T#1} }
\newcommand{\eT}{\end{Theorem}}

\newcommand{\bProposition}[1]{
\begin{Proposition} \label{P#1}}
\newcommand{\eP}{\end{Proposition}}

\newcommand{\bLemma}[1]{
\begin{Lemma} \label{L#1} }
\newcommand{\eL}{\end{Lemma}}

\newcommand{\bCorollary}[1]{
\begin{Corollary} \label{C#1} }
\newcommand{\eC}{\end{Corollary}}

\newcommand{\bFormula}[1]{
\begin{equation} \label{#1}}
\newcommand{\eF}{\end{equation}}

\newcommand{\Ov}[1]{\overline{#1}}
\newcommand{\DC}{C^\infty_c}

\newcommand{\vr}{\varrho}
\newcommand{\vre}{\vr_\ep}

\newcommand{\vue}{\vu_\ep}

\newcommand{\vu}{\vc{u}}
\newcommand{\vc}[1]{{\bf #1}}
\newcommand{\qed}{\bigskip \rightline {Q.E.D.} \bigskip}
\newcommand{\Div}{{\rm div}_x}
\newcommand{\Grad}{\nabla_x}

\newcommand{\tn}[1]{\mbox {\F #1}}
\newcommand{\dx}{{\rm d} {x}}
\newcommand{\vph}{{\boldsymbol \varphi}}
\newcommand{\dt}{{\rm d} t }

\newcommand{\dxdt}{\dx \ \dt}

\newcommand{\bProof}{{\bf Proof: }}

\newcommand{\vx}{\vc{x}}

\newcommand{\ep}{\varepsilon}

\font\F=msbm10 scaled 1000

\definecolor{grey}{rgb}{0.85,0.85,0.85}
\date{}

\long\def\greybox#1{%
    \newbox\contentbox%
    \newbox\bkgdbox%
    \setbox\contentbox\hbox to \hsize{%
        \vtop{
            \kern\columnsep
            \hbox to \hsize{%
                \kern\columnsep%
                \advance\hsize by -2\columnsep%
                \setlength{\textwidth}{\hsize}%
                \vbox{
                    \parskip=\baselineskip
                    \parindent=0bp
                    #1
                }%
                \kern\columnsep%
            }%
            \kern\columnsep%
        }%
    }%
    \setbox\bkgdbox\vbox{
        \color{grey}
        \hrule width  \wd\contentbox %
               height \ht\contentbox %
               depth  \dp\contentbox
        \color{black}
    }%
    \wd\bkgdbox=0bp%
    \vbox{\hbox to \hsize{\box\bkgdbox\box\contentbox}}%
    \vskip\baselineskip%
}

\begin{document}


\title{Weak solutions to the barotropic Navier-Stokes system with slip boundary conditions in time dependent domains}
\author{Eduard Feireisl \thanks{The work was supported by Grant 201/09/0917 of GA \v CR.} \and Ond{\v r}ej Kreml \thanks{The work was supported by Grant P201/11/1304 of GA \v CR.} \and {\v S}{\' a}rka Ne{\v c}asov{\' a}\thanks{The work was supported by Grant P201/11/1304 of GA \v CR.} \and Ji{\v r}{\' \i} \ Neustupa \thanks{The work was supported by Grant 201/08/0012 of GA \v CR.} \and Jan Stebel \thanks{The work was supported by Grant P201/11/1304 of GA \v CR.}}

\maketitle

\bigskip

\centerline{Institute of Mathematics of the Academy of Sciences of
the Czech Republic} \centerline{\v Zitn\' a 25, 115 67 Praha 1,
Czech Republic}

\medskip

\begin{abstract}

We consider the compressible (barotropic) Navier-Stokes system on time-dependent domains, supplemented with slip
boundary conditions. Our approach is based on penalization of the boundary behaviour, viscosity, and the pressure
in the weak formulation. Global-in-time weak solutions are obtained.

\end{abstract}

\medskip


\medskip

\section{Introduction}
\label{i}

Problems involving the motion of solid objects in fluids occur frequently in various applications of continuum fluid dynamics,
where the boundary conditions on the interfaces play a crucial role. Besides the commonly used \emph{no-slip} condition, where
the velocity of the fluid coincides with that of the adjacent solid body, various slip-like conditions have been proposed to handle the situations
in which the no-slip scenario fails to produce a correct description of the fluid boundary behavior, see Bul\' \i \v cek, M\' alek and Rajagopal  \cite{BMR1}, Priezjev a Troian \cite{PRTR} and the references therein. For \emph{viscous} fluids, Navier proposed the boundary conditions in the form
\bFormula{i1}
\left[ \tn{S} \vc{n} \right]_{\rm tan} + \kappa \left[ \vu - \vc{V} \right]_{\rm tan}|_{\Gamma_\tau} = 0, \ \kappa \geq 0,
\eF
where $\tn{S}$ is the viscous stress tensor, $\kappa$ represents a ``friction'' coefficient, $\vu$ and $\vc{V}$ denote the fluid and solid body velocities, respectively, and $\Gamma_\tau$ is the position of the interface at a time $\tau$, with the outer normal vector $\vc{n}$. If $\kappa = 0$, we obtain the \emph{complete slip}
while the asymptotic limit $\kappa \to \infty$ gives rise to the standard no-slip boundary conditions.

Besides their applications in ``thin'' domains occurring in nanotechnologies (see Qian, Wang and Sheng \cite{QIWA}), the slip boundary conditions are particularly relevant for dense viscous gases (see Coron \cite{Cor}), described by means of the standard \emph{Navier-Stokes system}:
\bFormula{i2}
\partial_t \vr + \Div (\vr \vu) = 0,
\eF
\bFormula{i3}
\partial_t (\vr \vu) + \Div (\vr \vu \otimes \vu) + \Grad p(\vr) = \Div \tn{S}(\Grad \vu) + \vr \vc{f},
\eF
where $\vr$ is the density, $p = p(\vr)$ the (barotropic) pressure, $\vc{f}$ a given external force, and $\tn{S}$ is determined by the standard
\emph{Newton rheological law}
\bFormula{i4}
\tn{S} (\Grad \vu) = \mu \left( \Grad \vu + \Grad^t \vu - \frac{2}{3} \Div \vu \tn{I} \right) + \eta \Div \vu \tn{I},\ \mu > 0,\ \eta \geq 0.
\eF

The boundary of the domain $\Omega_t$ occupied by the fluid is described by means of a \emph{given} velocity field $\vc{V}(t,\vx)$, where $t \geq 0$ and
$\vx \in R^3$. More specifically, assuming $\vc{V}$ is regular, we solve the associated system of differential equations
\bFormula{i5}
\frac{{\rm d}}{{\rm d}t} \vc{X}(t, \vx) = \vc{V} \Big( t, \vc{X}(t, \vx) \Big),\ t > 0,\ \vc{X}(0, \vx) = \vx,
\eF
and set
\[
\Omega_\tau = \vc{X} \left( \tau, \Omega_0 \right), \ \mbox{where} \ \Omega_0 \subset R^3 \ \mbox{is a given domain},\
\Gamma_\tau = \partial \Omega_\tau, \ \mbox{and}\ Q_\tau = \{ (t,x) \ |\ t \in (0,\tau), \ x \in \Omega_\tau \}.
\]

In addition to (\ref{i1}), we assume that the boundary $\Gamma_\tau$ is impermeable, meaning
\bFormula{i6}
(\vu - \vc{V}) \cdot \vc{n} |_{\Gamma_\tau} = 0 \ \mbox{for any}\ \tau \geq 0.
\eF

Finally, the problem (\ref{i1} - \ref{i6}) is supplemented by the initial conditions
\bFormula{i7}
\vr(0, \cdot) = \vr_0 ,\ (\vr \vu) (0, \cdot) = (\vr \vu)_0 \ \mbox{in}\ \Omega_0.
\eF

Our main goal is to show \emph{existence} of global-in-time \emph{weak solutions} to problem (\ref{i1} - \ref{i7}) for any finite energy initial data. The existence theory for the barotropic Navier-Stokes system on \emph{fixed} spatial domains in the framework of weak solutions was developed
in the seminal work by Lions \cite{LI4}, and later extended in \cite{FNP} to a class of physically relevant pressure-density state equations.
The investigation of \emph{incompressible} fluids in {time dependent} domains started with a seminal paper of Ladyzhenskaya \cite{LAD2}, see also 
\cite{Neust1}, \cite{NeuPen1}, \cite{NeuPen2} for more recent results in this direction.

\emph{Compressible} fluid flows in time dependent domains, supplemented with the \emph{no-slip} boundary conditions, were examined in \cite{FeNeSt} by means of Brinkman's penalization method. However, applying a penalization method to the slip boundary conditions is more delicate. Unlike for \emph{no-slip},
where the fluid velocity coincides with the field $\vc{V}$ outside $\Omega_\tau$, it is only its normal component $\vu \cdot \vc{n}$ that can be controlled in the case of \emph{slip}. In particular, given the rather poor {\it a priori} bounds available in the class of weak solutions,
we lose control over the boundary behavior of the normal stress $\tn{S} \vc{n}$ involved in Navier's condition (\ref{i1}).

A rather obvious penalty approach to slip conditions for \emph{stationary incompressible} fluids was proposed by Stokes and Carey \cite{StoCar}.
In the present setting, the variational (weak) formulation of the momentum equation is supplemented by a singular forcing term
\bFormula{i8}
\frac{1}{\ep} \int_0^T \int_{\Gamma_t}  (\vu - \vc{V} ) \cdot \vc{n} \ \vph \cdot \vc{n} \ {\rm dS}_{x} \ \dt,\ \ep > 0 \ \mbox{small},
\eF
penalizing the normal component of the velocity on the boundary of the fluid domain. In the time-dependent geometries, the penalization can be applied
in the \emph{interior} of a fixed reference domain, however, the resulting limit system consists of \emph{two} fluids separated by impermeable boundary and coupled through the tangential components of normal stresses. In such a way, an extra term is produced acting on the fluid by its ``complementary'' part outside $\Omega_\tau$. In order to eliminate these extra stresses, we use the following three level penalization scheme:
\begin{enumerate}
\item In addition to (\ref{i8}), we introduce a \emph{variable} shear viscosity coefficient $\mu = \mu_\omega$, where $\mu_\omega$ remains strictly positive in the fluid domain $Q_T$ but vanishes in the solid domain $Q_T^c$ as $\omega \to 0$.

\item Similarly to the existence theory developed in \cite{FNP}, we introduce the \emph{artificial pressure}
\[
p_\delta(\varrho) = p(\varrho) + \delta \varrho^\beta,\ \beta \geq 2,\ \delta > 0,
\]
in the momentum equation (\ref{i3}).
\item Keeping $\ep, \delta,\ \omega > 0$ fixed, we solve the modified problem in a (bounded) reference domain $B \subset R^3$ chosen in such a way that
\[
\Ov{ \Omega}_\tau \subset B \ \mbox{for any}\ \tau \geq 0.
\]
To this end, we adapt the existence theory for the compressible Navier-Stokes system with variable viscosity coefficients developed in \cite{EF71}.
\item We take the initial density $\varrho_0$ vanishing outside $\Omega_0$, and letting $\ep \to 0$ for fixed $\delta , \ \omega > 0$ we obtain a ``two-fluid'' system, where the density vanishes in the solid part $\left( (0,T) \times B \right) \setminus Q_T$ of the reference domain.

\item Letting the viscosity vanish in the solid part, we perform the limit $\omega \to 0$, where the extra stresses disappear in the limit system. The desired conclusion results from the final limit process $\delta \to 0$.

\end{enumerate}

The paper is organized as follows. In Section \ref{m}, we introduce all necessary preliminary material including a weak formulation of the problem and state the main result. Section \ref{p} is devoted to the penalized problem and to uniform bounds and existence of solutions at the starting level of approximations.
In Section \ref{s}, the singular limits for $\ep \to 0$, $\omega \to 0$, and $\delta \to 0$ are preformed successively. Section \ref{d} discusses possible extensions and applications of the method.

\section{Preliminaries, weak formulation, main result}
\label{m}

In the weak formulation, it is convenient that the equation of continuity (\ref{i2}) holds in the whole physical space $R^3$ provided the density
$\vr$ was extended to be zero outside the fluid domain, specifically
\bFormula{m1}
\int_{\Omega_\tau} \vr \varphi (\tau, \cdot) \ \dx - \int_{\Omega_0} \vr_0 \varphi (0, \cdot) \ \dx =
\int_0^\tau \int_{ \Omega_t} \left( \vr \partial_t \varphi + \vr \vu \cdot \Grad \varphi \right) \ \dxdt
\eF
for any $\tau \in [0,T]$ and any test function $\varphi \in \DC([0,T] \times R^3)$. Moreover, equation (\ref{i2}) is also satisfied in the sense of
renormalized solutions introduced by DiPerna and Lions \cite{DL}:
\bFormula{m2}
\int_{\Omega_\tau} b(\vr) \varphi (\tau, \cdot) \ \dx - \int_{\Omega_0} b(\vr_0) \varphi (0, \cdot) \ \dx =
\int_0^\tau \int_{ \Omega_t} \left( b(\vr) \partial_t \varphi + b(\vr) \vu \cdot \Grad \varphi +
\left( b(\vr)  - b'(\vr) \vr \right) \Div \vu \varphi \right) \ \dxdt
\eF
for any $\tau \in [0,T]$, any $\varphi \in \DC([0,T] \times R^3)$, and any $b \in C^1 [0, \infty)$, $b(0) = 0$, $b'(r) = 0$ for large $r$.
Of course, we suppose that $\vr \geq 0$ a.a. in $(0,T) \times R^3$.

Similarly, the momentum equation (\ref{i3}) is replaced by a family of integral identities
\bFormula{m3}
\int_{\Omega_\tau} \vr \vu \cdot \vph (\tau, \cdot) \ \dx - \int_{\Omega_0} (\vr \vu)_0 \cdot \vph (0, \cdot) \ \dx
\eF
\[
= \int_0^\tau \int_{\Omega_t} \left( \vr \vu \cdot \partial_t \vph + \vr [\vu \otimes \vu] : \Grad \vph + p(\vr) \Div \vph
- \tn{S} (\Grad \vu) : \Grad \vph + \vr \vc{f} \cdot \vph \right) \ \dxdt
\]
for any $\tau \in [0,T]$ and any test function $\vph \in \DC([0,T] \times R^3 ; R^3)$ satisfying
\bFormula{m4}
\vph \cdot \vc{n}|_{\Gamma_\tau} = 0 \ \mbox{for any} \ \tau \in [0,T].
\eF

Finally, the impermeability condition (\ref{i6}) is satisfied in the sense of traces, specifically,
\bFormula{m5}
\vu \in L^2(0,T; W^{1,2}(R^3; R^3)) \ \mbox{and}\ (\vu - \vc{V}) \cdot \vc{n}  (\tau , \cdot)|_{\Gamma_\tau}  = 0 \ \mbox{for a.a.}\ \tau \in [0,T].
\eF

At this stage, we are ready to state the main result of the present paper:

\greybox{

\bTheorem{m1}
Let $\Omega_0 \subset R^3$ be a bounded domain of class $C^{2 + \nu}$, and let $\vc{V} \in C^1([0,T]; C^{3}_c (R^3;R^3))$ be given.
Assume that the pressure $p \in C[0, \infty) \cap C^1(0, \infty)$ satisfies
\[
p(0) = 0,\ p'(\vr) > 0 \ \mbox{for any}\ \vr > 0,\ \lim_{\vr \to \infty} \frac{p'(\vr)}{\vr^{\gamma - 1}} = p_\infty > 0
\ \mbox{for a certain}\ \gamma > 3/2.
\]
Let the initial data satisfy
\[
\vr_0 \in L^\gamma (R^3),\ \vr_0 \geq 0, \ \vr_0 \not\equiv 0,\ \vr_0|_{R^3 \setminus \Omega_0} = 0,\
(\vr \vu)_0 = 0 \ \mbox{a.a. on the set} \ \{ \vr_0 = 0 \} ,\ \int_{\Omega_0} \frac{1}{\vr_0} |(\vr \vu)_0 |^2 \ \dx < \infty.
\]

Then the problem (\ref{i1} - \ref{i8}) admits a weak solution on any time interval $(0,T)$ in the sense specified through (\ref{m1} - \ref{m5}).
\eT

}

The rest of the paper is devoted to the proof of Theorem \ref{Tm1}.

\section{Penalization}
\label{p}

For the sake of simplicity,
we restrict ourselves to the case $\kappa = 0$, $\eta = 0$, and $\vc{f} = 0$. As we shall see, the main ideas of the proof presented below require only
straightforward modifications to accommodate the general case.

\subsection{Penalized problem - weak formulation}

Choosing $R > 0$ such that
\bFormula{p1}
\vc{V} |_{[0,T] \times \{ |\vx| > R \} } = 0 ,\ \ \Ov{ \Omega }_0 \subset \{ |\vx| < R \}
\eF
we take the reference domain $B = \{ |\vx| < 2R \}$.

Next, the shear viscosity coefficient $\mu = \mu_{\omega}(t,\vx)$ is taken such that
\bFormula{p2}
\mu_\omega \in \DC ([0,T] \times R^3),\ 0 < \underline{\mu}_\omega \leq \mu_\omega (t,x) \leq \mu \ \mbox{in}\ [0,T] \times B, \
\mu_{\omega}(\tau, \cdot)|_{ \Omega_\tau } = \mu \ \mbox{for any} \ \ \tau \in [0,T].
\eF

Finally, we define modified initial data so that
\bFormula{data1}
\vr_0 = \varrho_{0, \delta},\ \vr_{0, \delta} \geq 0, \ \vr_{0,\delta} \not\equiv 0,\ \vr_{0, \delta}|_{R^3 \setminus \Omega_0} = 0,\ \int_{B}
\left( \vr_{0, \delta}^\gamma +
\delta \vr_{0, \delta}^\beta \right) \dx \leq c,
\eF
\bFormula{data2}
(\vr \vu)_0 = (\vr \vu)_{0, \delta},\ (\vr \vu)_{0, \delta} = 0 \ \mbox{a.a. on the set} \ \{ \vr_{0, \delta} = 0 \} ,\ \int_{\Omega_0} \frac{1}{\vr_{0, \delta}} |(\vr \vu)_{0,\delta} |^2 \ \dx \leq c.
\eF

The weak formulation of the \emph{penalized problem} reads as follows:
\bFormula{p3}
\int_{B} \vr \varphi (\tau, \cdot) \ \dx - \int_{B} \vr_0 \varphi (0, \cdot) \ \dx =
\int_0^\tau \int_{B} \left( \vr \partial_t \varphi + \vr \vu \cdot \Grad \varphi \right) \ \dxdt
\eF
for any $\tau \in [0,T]$ and any test function $\varphi \in \DC([0,T] \times R^3)$;
\bFormula{p4}
\int_{B} \vr \vu \cdot \vph (\tau, \cdot) \ \dx - \int_{B} (\vr \vu)_0 \cdot \vph (0, \cdot) \ \dx
\eF
\[
= \int_0^\tau \int_{B} \left( \vr \vu \cdot \partial_t \vph + \vr [\vu \otimes \vu] : \Grad \vph + p(\vr) \Div \vph + \delta \vr^\beta \Div \vph
- \mu_\omega \left( \Grad \vu + \Grad^t \vu - \frac{2}{3} \Div \vu \tn{I} \right) : \Grad \vph \right) \ \dxdt
\]
\[
+ \frac{1}{\ep} \int_0^\tau \int_{ \Gamma_t } \left( (\vc{V} - \vu ) \cdot \vc{n} \ \vph \cdot \vc{n} \right) \ {\rm dS}_{x} \ \dt
\]
for any $\tau \in [0,T]$ and any test function $\vph \in \DC([0,T] \times B ; R^3)$,
where $\vu \in L^2(0,T; W^{1,2}_0 (B; R^3))$, meaning
$\vu$ satisfies the no-slip boundary condition
\bFormula{p5}
\vu|_{\partial B} = 0 \ \mbox{in the sense of traces}.
\eF
Here, $\ep$, $\delta$, and $\omega$ are positive parameters. The choice of the no-slip boundary condition (\ref{p5}) is not essential.

The \emph{existence} of global-in-time solutions to the penalized problem can be shown by means of the method developed in \cite{EF71} to handle the nonconstant viscosity coefficients. Indeed, for $\ep > 0$ fixed, the extra penalty term in (\ref{p4}) can be treated as a ``compact'' perturbation.
In addition, solutions can be constructed satisfying the \emph{energy inequality}
\bFormula{p6}
\int_B \left( \frac{1}{2} \vr |\vu|^2 + P(\vr) + \frac{\delta}{\beta - 1} \vr^\beta \right)(\tau, \cdot) \ \dx +
\frac{1}{2} \int_0^\tau \int_B \mu_\omega \left| \Grad \vu + \Grad^t \vu - \frac{2}{3} \Div \vu \tn{I} \right|^2 \ \dxdt
\eF
\[
+
\frac{1}{\ep} \int_0^\tau \int_{\Gamma_t} \left[ \left( \vu - \vc{V} \right) \cdot \vc{n} \right] \vu \cdot \vc{n} \ {\rm dS}_x \ \dt
\leq \int_B \left( \frac{1}{2 \vr_{0,\delta} } |(\vr \vu)_{0,\delta} |^2 + P(\vr_{0, \delta}) + \frac{\delta}{\beta - 1} \vr_{0,\delta}^\beta \right) \ \dx,
\]
where
\[
P(\vr) = \vr \int_1^\vr \frac{p(z)}{z^2} \ {\rm d}z .
\]
Note that the quantity on the right-hand side of (\ref{p6}) representing the total energy of the system is finite because of
(\ref{data1}), (\ref{data2}).

In addition, since $\beta \geq 2$, the density is square integrable and we may use the regularization technique of DiPerna and Lions \cite{DL}
to deduce the renormalized version of (\ref{p3}), namely
\bFormula{p6a}
\int_{B} b(\vr) \varphi (\tau, \cdot) \ \dx - \int_{B} b(\vr_0) \varphi (0, \cdot) \ \dx =
\int_0^\tau \int_{B} \left( b(\vr) \partial_t \varphi + b(\vr) \vu \cdot \Grad \varphi +
\left( b(\vr)  - b'(\vr) \vr \right) \Div \vu \varphi \right) \ \dxdt
\eF
for any $\varphi$ and $b$ as in (\ref{m2}).

\subsection{Modified energy inequality and uniform bounds}

Since the vector field $\vc{V}$ vanishes on the boundary $\partial B$, it may be used as a test function in (\ref{p4}).  Combining the resulting expression with the energy inequality
(\ref{p6}), we obtain
\bFormula{p7}
\int_B \left( \frac{1}{2} \vr |\vu|^2 + P(\vr) + \frac{\delta}{\beta - 1} \vr^\beta \right)(\tau, \cdot) \ \dx +
\frac{1}{2} \int_0^\tau \int_B \mu_\omega \left| \Grad \vu + \Grad^t \vu - \frac{2}{3} \Div \vu \tn{I} \right|^2 \ \dxdt
\eF
\[
+
\frac{1}{\ep} \int_0^\tau \int_{\Gamma_t} \left| \left( \vu - \vc{V} \right) \cdot \vc{n} \right|^2 \ {\rm dS}_x \ \dt
\leq \int_B \left( \frac{1}{2 \vr_{0,\delta}} |(\vr \vu)_{0,\delta} |^2 + P(\vr_{0, \delta}) + \frac{\delta}{\beta - 1} \vr_{0, \delta}^\beta \right) \ \dx
\]
\[
+ \int_B \Big(  (\vr \vu \cdot \vc{V}) (\tau, \cdot) - (\vr \vu)_{0,\delta} \cdot \vc{V}(0, \cdot) \Big)
\ \dx
\]
\[
+ \int_0^\tau \int_B \left( \mu_\omega \left(\Grad \vu + \Grad^t \vu - \frac{2}{3} \Div \vu \tn{I} \right) : \Grad \vc{V} - \vr \vu \cdot \partial_t \vc{V} - \vr \vu \otimes \vu : \Grad \vc{V} -
p(\vr) \Div \vc{V} - \frac{ \delta }{\beta - 1} \vr^\beta \Div \vc{V} \right)  \dxdt.
\]

Since the vector field $\vc{V}$ is regular and since
\[
p(\vr) \leq c(1 + P(\vr)) \ \mbox{for all} \ \vr \geq 0,
\]
relation (\ref{p7}) gives rise to the following bounds
\emph{independent} of the parameters $\ep$, $\delta$, and $\omega$:
\bFormula{p8}
{\rm ess} \sup_{t \in (0,T)} \| \sqrt{\vr} \vu (t , \cdot)  \|_{L^2 (B;R^3)} \leq c,
\eF
\bFormula{p9}
{\rm ess} \sup_{t \in (0,T)} \int_B P(\vr)(t , \cdot) \ \dx \leq c \ \mbox{yielding} \ {\rm ess} \sup_{t \in (0,T)} \left\| \vr (t, \cdot) \right\|_{L^\gamma(B)} \leq c,
\eF
\bFormula{p10}
{\rm ess} \sup_{t \in (0,T)} \delta \| \vr (t, \cdot) \|^\beta_{L^\beta(B)} \leq c,
\eF
\bFormula{p11}
\int_0^T \int_B \mu_{\omega} \left| \Grad \vu + \Grad^t \vu - \frac{2}{3} \Div \vu \tn{I} \right|^2 \ \dxdt \leq c,
\eF
and
\bFormula{p12}
\int_0^T \int_{\Gamma_t} \left| (\vu - \vc{V}) \cdot \vc{n}  \right|^2 \ {\rm dS}_{x} \ \dt  \leq \ep c.
\eF

Finally, we note that the total mass is conserved, meaning
\bFormula{p13}
\int_B \vr (\tau, \cdot) \ \dx = \int_B \vr_{0, \delta} \ \dx = \int_{\Omega_0} \vr_{0, \delta} \ \dx \leq c \ \mbox{for any}\ \tau \in [0,T].
\eF
Thus, relations (\ref{p8}), (\ref{p11}), (\ref{p13}), combined with the generalized version of Korn's inequality (see \cite[Theorem 10.17]{FeNo6}),
imply that
\bFormula{p14}
\int_0^T \| \vu (t, \cdot) \|^2_{W^{1,2}_0 (B;R^3)} \leq c(\omega).
\eF

\subsection{Pressure estimates}

Since the surfaces $\Gamma_\tau$ are determined {\it a priori}, we can use the technique based on the so-called Bogovskii operator to deduce the uniform bounds
\bFormula{p15}
\int\int_{\mathcal K} \Big( p(\vr) \vr^\nu + \delta \vr^{\beta + \nu} \Big) \ \dxdt \leq c(\mathcal {K}) \ \mbox{for a certain} \ \nu > 0
\eF
for any compact $\mathcal{K} \subset [0,T] \times \Ov{B}$ such that
\[
\mathcal{K} \cap \left( \cup_{ \tau \in [0,T] } \Big( \{ \tau \} \times \Gamma_\tau \Big) \right) = \emptyset,
\]
see \cite{FP13} for details.

Note that due to the fact that the boundaries $\Gamma_\tau$ change with time, \emph{uniform estimates} like (\ref{p15}) on the whole space time cylinder $(0,T) \times B$ seem to be a delicate matter. On the other hand, the mere \emph{equi-integrability} of the pressure could be shown by the method based on special test functions used in \cite{FeNeSt}.

\section{Singular limits}

\label{s}

In this section, we perform successively the singular limits $\ep \to 0$, $\omega \to 0$, and $\delta \to 0$.

\subsection{Penalization limit}

Keeping the parameters $\delta$, $\omega$ fixed, our goal is to let $\ep \to 0$ in (\ref{p3}), (\ref{p4}).
Let $\{ \vre, \vue \}$ be the corresponding weak solution of the perturbed problem constructed in the previous section. To begin, the estimates (\ref{p9}), (\ref{p14}),
combined with the equation of continuity (\ref{p3}), imply that
\[
\vre \to \vr \ \mbox{in}\ C_{\rm weak}([0,T] ; L^\gamma (B)),
\]
and
\[
\vue \to \vu \ \mbox{weakly in}\ L^2(0,T; W^{1,2}_0 (B, R^3))
\]
at least for suitable subsequences, where, as a direct consequence of (\ref{p12}),
\bFormula{s1}
( \vu - \vc{V} ) \cdot \vc{n} (\tau, \cdot) |_{\Gamma_\tau} = 0 \ \mbox{for a.a.}\ \tau \in [0,T].
\eF

Consequently, in accordance with (\ref{p8}), (\ref{p9}) and the compact embedding $L^\gamma(B) \hookrightarrow\hookrightarrow W^{-1,2}(B)$,
we obtain
\bFormula{s1a}
\vre \vue  \to \vr \vu \ \mbox{weakly-(*) in}\ L^\infty(0,T; L^{2 \gamma / (\gamma + 1)}(B;R^3)),
\eF
and, thanks to the embedding $W^{1,2}_0 (B) \hookrightarrow L^6(B)$,
\[
\vre \vue \otimes \vue \to \Ov{ \vr \vu \otimes \vu } \ \mbox{weakly in}\ L^2(0,T; L^{6 \gamma / (4 \gamma + 3)} (B; R^3)),
\]
where we have used the bar to denote a weak limit of a composed function.

Finally, it follows from the momentum equation (\ref{p4}) that
\[
\vre \vue \to \vr \vu \ \mbox{in} \ C_{\rm weak}([T_1, T_2]; L^{2 \gamma / (\gamma + 1)}(O;R^3))
\]
for any space-time cylinder
\[
(T_1, T_2) \times {O} \subset [0,T] \times B,\  [T_1, T_2] \times \Ov{O} \cap \cup_{\tau \in [0,T]} \left( \{ \tau \} \times \Gamma_\tau \right) =
\emptyset.
\]
Seeing that $L^{2 \gamma / (\gamma + 1)}(B) \hookrightarrow\hookrightarrow W^{-1,2}(B)$ we conclude, exactly as in (\ref{s1a}), that
\[
\Ov{\vr \vu \otimes \vu} = \vr \vu \otimes \vu \ \mbox{a.a. in}\ (0,T) \times B.
\]

Passing to the limit in (\ref{p3})  we obtain
\bFormula{s2}
\int_{B} \vr \varphi (\tau, \cdot) \ \dx - \int_{B} \vr_{0,\delta} \varphi (0, \cdot) \ \dx =
\int_0^\tau \int_{B} \left( \vr \partial_t \varphi + \vr \vu \cdot \Grad \varphi \right) \ \dxdt
\eF
for any $\tau \in [0,T]$ and any test function $\varphi \in \DC([0,T] \times R^3)$.

The limit in the momentum equation (\ref{p4}) is more delicate. Since we have at hand only the \emph{local estimates} (\ref{p15}) on the pressure,
we have to restrict ourselves to the class of test functions
\bFormula{s4}
\vph \in C^1([0,T]; W^{1, \infty}_0 (B; R^3)),\ {\rm supp}[ \Div \vph (\tau, \cdot)] \cap \Gamma_\tau = \emptyset,\
\vph \cdot \vc{n}|_{\Gamma_t} = 0 \ \mbox{for all}\ \tau \in [0,T].
\eF

Passing to the limit in (\ref{p4}), we obtain
\bFormula{s3}
\int_{B} \vr \vu \cdot \vph (\tau, \cdot) \ \dx - \int_{B} (\vr \vu)_{0, \delta} \cdot \vph (0, \cdot) \ \dx
\eF
\[
= \int_0^\tau \int_{B} \left( \vr \vu \cdot \partial_t \vph + \vr [\vu \otimes \vu] : \Grad \vph + \Ov{p(\vr)} \Div \vph + \delta \Ov{\vr^\beta} \Div \vph
- \mu_\omega \left( \Grad \vu + \Grad^t \vu - \frac{2}{3} \Div \vu \tn{I} \right) : \Grad \vph \right) \ \dxdt
\]
for any test function $\vph$ as (\ref{s4}). Note that the requirement of smoothness of $\vph$ postulated in (\ref{p4}) can be easily relaxed by means of a density argument.

Finally, we show pointwise (a.a.) convergence of the sequence $\{ \vre \}_{\ep > 0}$. To this end, we adopt the method developed in \cite{EF71}
to accommodate the variable viscosity coefficient $\mu_\omega$. The crucial observation is the effective viscous pressure identity that can be established exactly as
in \cite{EF71}:
\bFormula{s5}
\Ov{ p_\delta (\vr) T_k (\vr) } - \Ov{ p_\delta (\vr) } \ \Ov{ T_k (\vr) } = \frac{4}{3} \mu_\omega \left(
\Ov{ T_k (\vr) \Div \vu } - \Ov{T_k (\vr) } \Div \vu \right)
\eF
where we have denoted
\[
p_\delta (\vr) = p(\vr) + \delta \vr^\beta ,\ T_k(\vr) = \min \{ \vr , k \}.
\]
Similarly to the pressure estimates (\ref{p15}), identity (\ref{s5}) holds only on compact sets $\mathcal{K} \subset [0,T] \times B$
satisfying
\[
\mathcal{K} \cap \left( \cup_{ \tau \in [0,T] } \Big( \{ \tau \} \times \Gamma_\tau \Big) \right) = \emptyset.
\]
We recall that this step leans on the satisfaction of the renormalized equation (\ref{p6a}) for both $\vre$ and the limit $\vr$ that can be shown by the regularization procedure of DiPerna and Lions \cite{DL}.

Following \cite{EF71}, we introduce the \emph{oscillations defect measure}
\[
{\bf osc}_q [ \vre \to \vr] (\mathcal{K}) = \sup_{k \geq 0} \left( \limsup_{\ep \to 0}
\int_{\mathcal{K}} | T_k (\vre) - T_k(\vr) |^q \ \dxdt \right),
\]
and use (\ref{s5}) to conclude that
\bFormula{s6}
{\bf osc}_{\gamma + 1} [ \vre \to \vr] (\mathcal{K}) \leq c(\omega) < \infty,
\eF
where the constant $c$ is \emph{independent} of $\mathcal{K}$. Thus
\bFormula{s7}
{\bf osc}_{\gamma + 1} [ \vre \to \vr] ([0,T] \times B) \leq c(\omega),
\eF
which implies, by virtue of the procedure developed in \cite{EF70}, the desired conclusion
\bFormula{s8}
\vre \to \vr \ \mbox{a.a. in}\ (0,T) \times B.
\eF

In accordance with (\ref{s8}), the momentum equation (\ref{s3}) reads
\bFormula{s9}
\int_{B} \vr \vu \cdot \vph (\tau, \cdot) \ \dx - \int_{B} (\vr \vu)_{0,\delta} \cdot \vph (0, \cdot) \ \dx
\eF
\[
= \int_0^\tau \int_{B} \left( \vr \vu \cdot \partial_t \vph + \vr [\vu \otimes \vu] : \Grad \vph + {p(\vr)} \Div \vph + \delta {\vr^\beta} \Div \vph
- \mu_\omega \left( \Grad \vu + \Grad^t \vu - \frac{2}{3} \Div \vu \tn{I} \right) : \Grad \vph \right) \ \dxdt
\]
for any test function $\vph$ as in (\ref{s4}). In addition, as already observed, the limit solution $\{ \vr, \vu \}$ satisfies also the renormalized equation (\ref{p6a}).

\subsubsection{Fundamental lemma}
\label{a}

Our next goal is to use the specific choice of the initial data $\vr_{0, \delta}$ to get rid of the density-dependent terms in (\ref{s9}) supported by the ``solid'' part $\left( (0,T) \times B \right) \setminus Q_T$. To this end, we show the following result, rather obvious for regular solutions but a bit more delicate in the weak framework, that may be of independent interest.

\greybox{

\bLemma{a1}
Let $\vr \in L^\infty (0,T; L^2(B))$, $\vr \geq 0$,  $\vu \in L^2(0,T; W^{1,2}_0(B;R^3))$ be a weak solution of the equation of continuity,
specifically,
\bFormula{a1}
\int_{B} \Big( \vr (\tau, \cdot) \varphi (\tau, \cdot) - \vr_0 \varphi(0, \cdot) \Big) \dx
= \int_0^\tau \int_{B} \Big( \vr \partial_t \varphi + \vr \vu \cdot \Grad \varphi \Big) \ \dx \ \dt
\eF
for any $\tau \in [0,T]$ and any test function $\varphi \in C^1_c ([0,T] \times R^3)$.

In addition, assume that
\bFormula{a2}
( \vu - \vc{V} )(\tau, \cdot) \cdot \vc{n} |_{\Gamma_\tau}  = 0 \ \mbox{for a.a.}\ \tau \in (0,T),
\eF
and that
\[
\vr_0 \in L^2 (R^3), \ \vr_0 \geq 0,  \ \vr_0 |_{B \setminus \Omega_0} = 0.
\]

Then
\[
\vr(\tau, \cdot) |_{B \setminus \Omega_\tau} = 0 \ \mbox{for any}\ \tau \in [0,T].
\]

\eL

}

\bigskip

\bProof

We use the level set approach to describe the interface $\Gamma_\tau$, see Osher and Fedkiw \cite{OshFed}. To this end,
we introduce a function $d = d(t,x)$ defined as the unique solution of the transport equation
\[
\partial_t d + \Grad d(t,x) \cdot \vc{V} = 0 ,\ t > 0, \ x \in R^3,
\]
with the initial data
\[
d(0, x) = d_0(x) \in C^\infty (R^3), \ d_0(x) = \left\{ \begin{array}{l} > 0 \ \mbox{for}\ x \in B \setminus \Omega_0 ,\\ \\
                                                         < 0 \ \mbox{for}\ x \in \Omega_0 \cup (R^3 \setminus \Ov{B}) \end{array} \right. ,
                                                         \ \Grad d_0 \ne 0 \ \mbox{on} \ \Gamma_0.
\]
Note that the interface $\Gamma_\tau$ can be identified with a component of the level set $\{ d(\tau, \cdot) = 0 \}$, while the sets $B \setminus \Omega_\tau$
correspond to $\{ d (\tau, \cdot) > 0 \}$. Finally,
\bFormula{a3}
\Grad d(\tau,x) = \lambda (\tau,x) \vc{n}(x) ,\ \mbox{for any}\ x \in \Gamma_\tau,\ \lambda(\tau, x) \geq 0 \ \mbox{for}\ \tau \in [0,T].
\eF

For a given $\xi > 0$, we take
\[
\varphi = \left[ \min \left\{ \frac{1}{\xi} d ; 1 \right\} \right]^+
\]
as a (Lipschitz) test function in the variational formulation (\ref{a1}) to obtain
\bFormula{a3a}
\int_{B \setminus \Omega_\tau} \vr\varphi (\tau, \cdot) \ \dx =
\frac{1}{\xi} \int_0^{\tau} \int_{ \{  0 \leq d(t,x) < \xi \} } \Big( \vr \partial_t d + \vr \vu \cdot \Grad d \Big) \ \dx \ \dt.
\eF

Now, we have
\[
\vr \partial_t d + \vr \vu \cdot \Grad d = \vr \Big( \partial_t d + \vu \cdot \Grad d \Big) = \vr \left( \vc{V} - \vu \right) \cdot
\Grad d
\]
where, by virtue of hypothesis (\ref{a2}) and relation (\ref{a3}),
\bFormula{a4}
\left( \vc{V} - \vu \right) \cdot
\Grad d \in W^{1,2}_0 (B \setminus \Omega_t) \ \mbox{for a.a.}\ t \in (0, \tau).
\eF
Introducing
\[
\delta (t,x) = {\rm dist}_{R^3} [x , \partial ( B \setminus \Omega_t) ] \ \mbox{for}\ t \in [0, \tau],\ x \in B \setminus \Omega_t,
\]
we deduce from (\ref{a4}) and Hardy's inequality that
\[
\frac{1}{\delta} \left( \vc{V} - \vu \right) \cdot
\Grad d \in L^2 ([0,\tau] \times B \setminus Q_\tau).
\]

Finally, since $\vc{V}$ is regular, we have
\[
\frac{\delta(t,x) }{\xi} \leq c  \ \mbox{whenever}  \ 0 \leq d(t,x) < \xi;
\]
whence, letting $\xi \to 0$ in (\ref{a3a}), we obtain the desired conclusion
\[
\int_{B \setminus \Omega_\tau} \vr (\tau, \cdot) \ \dx =  0,
\]
where we have used the fact that $\vr \in L^\infty(0,T; L^2(B))$.

\qed

Thus, by virtue of Lemma \ref{La1}, the momentum equation (\ref{s9}) reduces to
\bFormula{s20}
\int_{\Omega_\tau } \vr \vu \cdot \vph (\tau, \cdot) \ \dx - \int_{\Omega_0} (\vr \vu)_{0, \delta} \cdot \vph (0, \cdot) \ \dx
\eF
\[
= \int_0^\tau \int_{\Omega_t} \left( \vr \vu \cdot \partial_t \vph + \vr [\vu \otimes \vu] : \Grad \vph + {p(\vr)} \Div \vph + \delta {\vr^\beta} \Div \vph
- \mu \left( \Grad \vu + \Grad^t \vu - \frac{2}{3} \Div \vu \tn{I} \right) : \Grad \vph \right) \ \dxdt
\]
\[
- \int_0^\tau \int_{ B \setminus \Omega_t } \mu_\omega \left( \Grad \vu + \Grad^t \vu - \frac{2}{3} \Div \vu \tn{I} \right) : \Grad \vph  \ \dxdt
\]
for any test function $\vph$ as in (\ref{s4}). We remark that it was exactly this step when we needed the extra pressure term $\delta \vr^\beta$
ensuring the density $\vr$ to be square integrable.

\subsection{Vanishing viscosity limit}

In order to get rid of the last integral in (\ref{s20}), we take the viscosity coefficient
\[
\mu_\omega = \left\{ \begin{array}{l} \mu = {\rm const} > 0 \ \mbox{in} \ Q_T , \\ \\
\mu_\omega \to 0 \ \mbox{a.a. in} \ ((0,T) \times B) \setminus Q_T. \end{array} \right.
\]
Denoting $\{ \vr_\omega, \vu_\omega\}$ the corresponding solution constructed in the previous section, we
may use (\ref{p11}) to deduce that
\bFormula{s22}
\int_0^T \int_{\Omega_t} \left| \Grad \vu_\omega + \Grad^t \vu_\omega - \frac{2}{3} \Div \vu_\omega \tn{I} \right|^2 \ \dxdt < c,
\eF
while
\[
\int_0^T \int_{ B \setminus \Omega_t } \mu_\omega \left| \Grad \vu_\omega + \Grad^t \vu_\omega - \frac{2}{3} \Div \vu_\omega \tn{I} \right|^2 \ \dxdt \leq c,
\]
where the latter estimates yields
\[
\int_0^\tau \int_{ B \setminus \Omega_t } \mu_\omega \left( \Grad \vu_\omega + \Grad^t \vu_\omega - \frac{2}{3} \Div \vu_\omega \tn{I} \right) : \Grad \vph  \ \dxdt =
\]
\[
\int_0^\tau \int_{ B \setminus \Omega_t } \sqrt{\mu_\omega} \sqrt{\mu_\omega} \left( \Grad \vu_\omega + \Grad^t \vu_\omega - \frac{2}{3} \Div \vu_\omega \tn{I} \right)
: \Grad \vph \ \dxdt \to 0 \ \mbox{as}\ \omega \to 0
\]
for any fixed $\vph$.

On the other hand, as we know from Lemma \ref{La1} that the density $\varrho_\omega$ is supported by the ``fluid'' region $Q_T$, we can still use (\ref{p8}), (\ref{s22}), together with Korn's inequality to obtain
\[
\int_0^T \int_{\Omega_t} | \Grad \vu_\omega |^2 \ \dxdt \leq c.
\]

Repeating step by step the arguments of the preceding section, we let $\omega \to 0$ to obtain the momentum equation in the form
\bFormula{s23}
\int_{\Omega_\tau} \vr \vu \cdot \vph (\tau, \cdot) \ \dx -  \int_{\Omega_0} (\vr \vu)_{0, \delta} \vph (0, \cdot) \ \dx
\eF
\[
= \int_0^\tau \int_{\Omega_t } \Big( \vr \vu \cdot \partial_t \vph + \vr [\vu \otimes \vu] : \Grad \vph  + p(\vr) \Div \vph
+ \delta \vr^\beta \Div \vph - \tn{S}(\Grad \vu) : \Grad \vph \Big) \ \dxdt
\]
for  any test function $\vph$ as in (\ref{s4}). Note that compactness of the density is now necessary only in the ``fluid'' part $Q_T$ so
a possible loss of regularity of $\vu_\omega$ outside $Q_T$ is irrelevant.

\subsection{Vanishing artificial pressure}

The final step is standard, we let $\delta \to 0$ in (\ref{s23}) to get rid of the artificial pressure term
$\delta \vr^\beta$ and to adjust the initial conditions, see \cite[Chapter 6]{EF70}. However, the momentum equation identity (\ref{s23}) holds only for the class of functions specified in (\ref{s4}). The last step of the proof of Theorem \ref{Tm1} is therefore to show that the class of admissible test functions can be extended by density arguments. This will be shown in the following part.

\subsubsection{Extending the class of test functions}

\label{ex}

Consider a test function $\vph \in \DC ([0,T] \times R^3 ; R^3)$ such that
\bFormula{e1}
\vph (\tau , \cdot) \cdot \vc{n} |_{\Gamma_\tau} = 0 \ \mbox{for any}\ \tau.
\eF
Our goal is to show the existence of an approximating sequence of functions $\vph_n$ belonging to the class specified in (\ref{s4})
and such that
\bFormula{e2}
\| \vph_n \|_{W^{1,\infty}((0,T) \times B;R^3)} \leq c,\ \vph_n \to \vph ,\ \partial_t \vph_n \to \partial_t \vph,
\ \mbox{and}\ \Grad \vph_n \to \Grad \vph \ \mbox{a.a. in} \ Q_T.
\eF
Combining (\ref{e2}) with Lebesgue dominated convergence theorem we may infer that $\vph$ belongs to the class of admissible test functions for
(\ref{m3}).

In other words, we have to find a suitable \emph{solenoidal extension} of the tangent vector field $\vph|_{\Gamma_\tau}$ inside
$\Omega_\tau$. Since $\Gamma_\tau$ is regular, there is an open neighborhood $\mathcal{U}_\tau$ of $\Gamma_\tau$ such that each point
$\vc{x} \in \mathcal{U}_\tau$ admits a \emph{single} closest point $\vc{b}_\tau(\vx) \in \Gamma_\tau$. We set
\[
\vc{h}(\tau, \vx) = \vph (\tau , \vc{b}_\tau(\vx)) \ \mbox{for all} \ \vx \in \mathcal{U}_\tau.
\]
Finally, we define
\[
\vc{w}(\tau , \vx) = \vc{h}(\tau, \vx) + \vc{g} (\tau, \vx),
\]
where
\[
\vc{g}(\tau, \vx) = 0 \ \mbox{whenever} \ \vx \in \Gamma_\tau,
\]
and, taking the local coordinate system at $\vx$ so that $e_3$ coincides with $\vx - \vc{b}_\tau(\vx)$, we set
\[
\vc{g}(\tau, \vx) = [0,0, g^3(\tau,\vx)], \partial_{x_3} g^3(\tau, \vx) = - \partial_{x_1} h^1(\tau, \vx) - \partial_{x_2} h^2 (\tau, \vx).
\]

We check that
\[
\Div \vc{w}(\tau, \cdot) = 0 \ \mbox{in}\ \mathcal{U}_\tau ,\ \vc{w}(\tau, \cdot)|_{\Gamma_\tau} = \vph (\tau , \cdot)|_{\Gamma_\tau}.
\]
Furthermore, extending $\vc{w}(\tau , \cdot)$ inside $\Omega_\tau$, we may use smoothness of $\vph$ and $\Gamma_\tau$ to conclude that
\[
\vc{w} \in W^{1,\infty} (Q_T).
\]

As a matter of fact, a (smooth) extension of $\vph$, \emph{solenoidal} in the whole domain $Q_T$, was constructed by
Shifrin \cite[Theorem 4]{Shif}.

Writing
\[
\vph = (\vph - \vc{w}) + \vc{w},
\]
we check that the field $\vc{w}$ belongs to the class (\ref{s4}), while
\[
(\vph - \vc{w})(\tau, \cdot)|_{\partial \Omega_\tau} = 0 \ \mbox{for any}\ \tau \geq 0.
\]
Thus, finally, it is a matter of routine to construct a sequence $\vc{a}_n$ such that
\[
\vc{a}_n \in \DC ([0,T] \times B; R^3),\ {\rm supp} [ \vc{a}_n (\tau, \cdot) ] \subset \Omega_\tau \ \mbox{for any}\ \tau \in [0,T],
\]
in particular $\vc{a}_n$ belongs to the class (\ref{s4}), and
\[
\| \vc{a}_n \|_{W^{1,\infty}((0,T) \times B;R^3)} \leq c,\ \vc{a}_n \to (\vph - \vc{w})  ,\ \partial_t \vc{a}_n \to \partial_t (\vph - \vc{w}),
\ \mbox{and}\ \Grad \vc{a}_n \to \Grad (\vph - \vc{w}) \ \mbox{a.a. in} \ Q_T.
\]
Clearly, the sequence
\[
\vph_n = \vc{a}_n + \vc{w}
\]
complies with (\ref{e2}).

We have completed the proof of Theorem \ref{Tm1}.

\section{Discussion}
\label{d}

The assumption on monotonicity of the pressure is not necessary, the same result can be obtained for a non-monotone pressure adopting the method
developed in \cite{EF61}.

As already pointed out, the technicalities of Section \ref{ex} could be avoided by means of the construction of special test functions
used in \cite{FeNeSt}. However,
it would be interesting to show that the pressure is bounded in some $L^q(Q_T)$, with $q > 1$, meaning that the estimate (\ref{p15}) holds in
$Q_T$.

As pointed out in the introduction, the general Navier slip conditions (\ref{i1}) are obtained introducing another boundary integral in the weak formulation, namely
\[
\int_0^T \int_{\Gamma_t}  \kappa (\vu - \vc{V}) \cdot \vph \ {\rm dS}_x \ \dt
\]
Taking $\kappa = \kappa(\vx)$ as a singular parameter, we can deduce results for mixed type no-slip - (partial) slip
boundary conditions prescribed on various components of $\Gamma_t$.

Last but not least, the method can be extended to unbounded (exterior) domains with prescribed boundary conditions
``at infinity''.

\def\ocirc#1{\ifmmode\setbox0=\hbox{$#1$}\dimen0=\ht0 \advance\dimen0
  by1pt\rlap{\hbox to\wd0{\hss\raise\dimen0
  \hbox{\hskip.2em$\scriptscriptstyle\circ$}\hss}}#1\else {\accent"17 #1}\fi}

\end{document}